\documentclass[11pt,twoside]{article}
\usepackage{amsfonts}
\usepackage{mathrsfs}
\usepackage{latexsym}
\usepackage{graphicx}
\usepackage{epsfig}

\normalsize \makeatletter

\begin{document}

\title{THE EXTENSION FOR MEAN CURVATURE FLOW WITH FINITE INTEGRAL CURVATURE IN RIEMANNIAN MANIFOLDS \footnote{
\newline\indent Research supported by the NSFC, Grant No. 10771187; the Trans-Century Training Programme Foundation for
Talents by the Ministry of Education of China; and the Natural
Science Foundation of Zhejiang Province, Grant No. 101037.
\newline \indent Keywords: mean curvature flow, Riemannian manifold, maximal existence
time, integral curvature.}}
\author{H{\footnotesize ONG-WEI} X{\footnotesize U},\ \ F{\footnotesize EI} Y{\footnotesize E}\ {\footnotesize AND}\ E{\footnotesize N-TAO} Z{\footnotesize HAO} \\
}
\date{}
\maketitle

%-----------------------------------------------------------------------------------------------------------------

\begin{abstract}
%\textbf{Abstract}:
We investigate the integral conditions to extend the mean curvature
flow in a Riemannian manifold. We prove that the mean curvature flow
solution with finite total mean curvature on a finite time interval
$[0,T)$ can be extended over time $T$. Moreover, we show that the
condition is optimal in some sense.

\end{abstract}

%-----------------------------------------------------------------------------------------------------------------

%-----------------------------------------------------------------------------------------------------------------

\section{Introduction}

\hspace*{6mm}Let $(M,g)$ be a compact $n$-dimensional manifold
without boundary, and let $F_t :M^n\rightarrow N^{n+1}$ be a
one-parameter family of smooth hypersurfaces immersed in a
Riemannian manifold $(N,h)$. We say that $M_t = F_t(M)$ is a
solution of the mean curvature flow if $F_t$ satisfies
\[
\left\{
\begin{array}{ccc}
\frac{\partial}{\partial t}F(x,t)&=&-H(x,t)\nu(x,t)\\
F(x,0)&=&F_0(x),
\end{array}\right.\]
where $F(x,t)=F_t(x)$, $H(x,t)$ is the mean curvature, $\nu(x,t)$ is
the unit outward normal vector, and $F_{0}$ is some given initial
hypersurface.

When the ambient space is the Euclidean space $\mathbb{R}^{n+1}$, G.
Huisken \cite{6} showed that the solution of the mean curvature flow
converges to a round point in a finite time for convex initial
hypersurface. He also proved that if the second fundamental form is
uniformly bounded, then the mean curvature flow can be extended over
time. If the ambient space is a Riemannian manifold, G. Huisken
\cite{7} proved the similar convergence theorem for certain initial
compact hypersurface and gave an sufficient condition to assure the
extension over time for mean curvature flow. Distinct from the above
pointwise conditions, in our previous work \cite{14} we investigated
the integral conditions to extend the mean curvature flow on closed
hypersurfaces in $\mathbb{R}^{n+1}$, which is optimal in some sense.
Almost at the same time, N. Le and N. \v{S}e\v{s}um \cite{8} studied
the same question independently with a different method.

In this paper, we study the mean curvature flow of hypersurfaces in
a Riemannian manifold with bounded geometry, which generalizes our
results in \cite{14}. We recall that a Riemannian manifold is said
to have bounded geometry if (i): the sectional curvature is bounded;
(ii): the first covariant derivative of the curvature tensor is
bounded; (iii): the injective radius is bounded from below by a
positive constant. In this paper we always assume that the ambient
space $N^{n+1}$ is a complete Riemannian manifold with bounded
geometry. We prove that when the space-time integration of the mean
curvature is finite and the second fundamental form is bounded from
below, the mean curvature flow can be extended.\\\\
\textbf{Theorem 1.1.} \emph{Let $F_t:\ M^n\longrightarrow N^{n+1}$
$(n\geq 3)$ be a solution of the mean curvature flow of closed
hypersurfaces on
a finite time interval $[0,T)$. If}\\
(1) \emph{there is a positive constant $C$ such that $h_{ij}\geq -C$
for $(x,t)\in M\times [0,T)$},\\
(2) \emph{$|| H||_{\alpha,M\times [0,T)}=\left(\int^T_0\int_M | H
|^\alpha d\mu dt\right)^{\frac{1}{\alpha}}<+\infty$ for some
$\alpha\geq n+2$,\\
then this flow can be extended over time $T$.}\\

Suppose that the sectional curvature $K_N$ of $N^{n+1}$ satisfies
$$-K_1\leq K_N\leq K_2,$$
where $K_1$ and $K_2$ are nonnegative constants. We also will prove the following theorem.\\\\
\textbf{Theorem 1.2.} \emph{Let $F_t:\ M^n\longrightarrow N^{n+1}$
$(n\geq 3)$ be a solution of the mean curvature flow of closed
hypersurfaces on
a finite time interval $[0,T)$. If}\\
(1) \emph{$H^2>n^2K_1$ at $t=0$},\\
(2) \emph{$|| H||_{\alpha,M\times [0,T)}=\left(\int^T_0\int_M | H
|^\alpha d\mu dt\right)^{\frac{1}{\alpha}}<+\infty$ for some
$\alpha\geq n+2$,\\
then this flow can be extended over time $T$.}\\

The following example shows that the condition $\alpha\geq n+2$ in
Theorem 1.1 and 1.2 is optimal when the ambient space is a complete simply connected space form.\\\\
\textbf{Example.} (i) For the case where $N^{n+1}=\mathbb{R}^{n+1}$,
set $\mathbb{S}^n=\{x\in\mathbb{R}^{n+1}:\sum^{n}_{i=1}x_i^2=1\}$.
Let $F$ be the standard isometric embedding of $\mathbb{S}^n$ into
$\mathbb{R}^{n+1}$. It is clear that $F(t)=\sqrt{1-2nt}F$ is the
solution to the mean curvature flow, where $T=\frac{1}{2n}$ is the
maximal existence time. By a simple computation, we have
$g_{ij}(t)=(1-2nt)g_{ij}$, $H(t)=\frac{n}{\sqrt{1-2nt}}$ and
$h_{ij}(t)>0$. Hence
\begin{eqnarray*}
|| H||_{\alpha,M\times [0,T)}&=&\left(\int^T_0\int_M
| H |^\alpha d\mu dt\right)^{\frac{1}{\alpha}}\\
&=&C_1\left(\int^T_0
(T-t)^{\frac{n-\alpha}{2}}dt\right)^{\frac{1}{\alpha}},
\end{eqnarray*}
where $C_1$ is a positive constant. It follows that

\[|| H||_{\alpha,M\times
[0,T)}\left\{ \begin{array}{ll}
=\infty,&\ for\ \ \alpha \geq n+2,\\
<\infty,&\ for\ \ \alpha < n+2.
\end{array}
\right.\] This implies that the condition $\alpha\geq n+2$ in
Theorem 1.1 and 1.2 is optimal when $N^{n+1}=\mathbb{R}^{n+1}$.\\\\
(ii) Let $\mathbb{F}^{n+1}(c)$ be a complete simply connected space
form with constant curvature $c$. We consider the case
$N^{n+1}=\mathbb{F}^{n+1}(c)$, where $c=\pm1$, that is,
$N^{n+1}=\mathbb{S}^{n+1}$ or $\mathbb{H}^{n+1}$. Let
$M=\mathbb{S}^n(r_0)$ be a total umbilical sphere of radius $r_0$ in
$N^{n+1}$ with constant mean curvature $H_0$ satisfying $H_0>0$ when
$c=1$, and $H_0^2>n^2$ when $c=-1$. Put
$d=\frac{H_0^2}{H^2_0+n^2c}$. Let $\mathbb{S}^{n}(r(t))$ be a sphere
with radius $r(t)=\frac{\sqrt{H_0^2+n^2c}}{\sqrt{H^2(t)+n^2c}}r_0$,
where $H(t)=\sqrt{\frac{n^2cde^{2nct}}{1-de^{2nct}}}$. Then
$\mathbb{S}^{n}(r(t))$ is a family of total umbilical spheres with
constant mean curvature $H(t)$, which satisfies the mean curvature
flow with initial value $M=\mathbb{S}^n(r_0)\subset N^{n+1}$. It is
clear that the maximal existence time is $T=-\frac{\ln d}{2n}$ when
$c=1$, and $T=\frac{\ln d}{2n}$ when $c=-1$, the second fundamental
form $h_{ij}$ satisfies $h_{ij}>0$, and the volume of
$\mathbb{S}^n(r(t))$ is
$V(t)=\bigg(\frac{H_0^2+n^2c}{H^2(t)+n^2c}\bigg)^{\frac{n}{2}}V_0$,
where $V_0$ is the volume of $\mathbb{S}^n(r_0)$. Hence
\begin{eqnarray*}||H||_{\alpha,M\times [0,T)}&=&\left(\int^T_0\int_M | H
|^\alpha d\mu dt\right)^{\frac{1}{\alpha}}\\
&=&\left(\int_0^T
H^{\alpha}(t)V(t)dt\right)^{\frac{1}{\alpha}}\\
&=&C_2\left(\int_0^T
({n^{2-n}de^{2nct}})^{\frac{\alpha}{2}}\bigg(\frac{1-de^{2nct}}{c}\bigg)^{\frac{n-\alpha}{2}}dt\right)^{\frac{1}{\alpha}},
\end{eqnarray*} where $C_2$ is a positive constant.
Since $({n^{2-n}de^{2nct}})^{\frac{\alpha}{2}}$ has positive upper
and lower bounds because of the finiteness of $T$, and the integral
\[\int_0^T\bigg(\frac{1-de^{2nct}}{c}\bigg)^{\frac{n-\alpha}{2}}dt\left\{ \begin{array}{ll}
=\infty,&\ for\ \ \alpha \geq n+2,\\
<\infty,&\ for\ \ \alpha < n+2,
\end{array}
\right.\] it follows that
\[|| H||_{\alpha,M\times
[0,T)}\left\{ \begin{array}{ll}
=\infty,&\ for\ \ \alpha \geq n+2,\\
<\infty,&\ for\ \ \alpha < n+2.
\end{array}
\right.\] This implies that the condition $\alpha\geq n+2$ in
Theorem 1.1 and 1.2 is optimal when $N^{n+1}=\mathbb{S}^{n+1}$ or
$\mathbb{H}^{n+1}$.

\section{Preliminaries}

\hspace*{6mm}Let $F_t :M^n\rightarrow N^{n+1}$ be a one-parameter
family of smooth hypersurfaces immersed in a Riemannian manifold
$N$. Denote by $g=\{g_{ij}\}$ and $A=\{h_{ij}\}$ the induced metric
and the second fundamental form of $M$ respectively, and $H$ is the
mean curvature of $M$, which is the trace of $A$. We put
$\bar{\nabla}$ and $\bar{R}ic$ be the connection and the Ricci
tensor of $N$, and $R_{ABCD}$, $A$, $B$, $C$, $D=0,1,\cdots ,n$, be
components of the curvature tensor of $N$ with respect to some local
coordinates such that $e_0=\nu$.

Firstly, we recall
some evolution equations (see \cite{3}, \cite{7} or \cite{16}).\\\\
\textbf{Lemma 2.1.} \emph{Along the mean curvature flow in
Riemannian manifolds, we have the following evolution equations}
\begin{eqnarray*}
\frac{\partial}{\partial t}g_{ij}&=&-2Hh_{ij},\\
\frac{\partial |A|^2}{\partial t}&=&\triangle|A|^2-2|\nabla A|^2+2|A|^2(|A|^2+\bar{R}ic(\nu,\nu))\\
&&-4(h^{ij}h^m_j\bar{R}_{mli}\
^{l}-h^{ij}h^{lm}\bar{R}_{milj})-2h^{ij}(\bar{\nabla}_j\bar{R}_{0li}\
^{l}+\bar{\nabla}_l\bar{R}_{0ij}\ ^{l}),\\
\frac{\partial}{\partial t}H&=&\triangle
H+H(|A|^2+\bar{R}ic(\nu,\nu)).
\end{eqnarray*}

We denote by $\omega_n$  the volume of the
unit ball in $\mathbb{R}^n$. The following Sobolev inequality can be found in \cite{9}.\\
\\
\textbf{Lemma 2.2.} \emph{Let $M^n\subset N^{n+p}$ be an
$n(\geq2)$-dimensional closed submanifold in a Riemannian manifold
$N^{n+p}$ with codimension $p\geq1$. Denote by $i_N$ the positive
lower bound of the injective radius of $N$ restricted on $M$. Assume
$K_N\leq b^2$ and let $h$ be a non-negative $C^1$ function on $M$.
Then}
$$\left(\int_Mh^{\frac{n}{n-1}}d\mu\right)^{\frac{n-1}{n}}\leq
C(n,\alpha)\int_M\Big[|\nabla h|+h|H|\Big]d\mu,$$ \emph{provided}
$$b^2(1-\alpha)^{-\frac{2}{n}}(\omega_n^{-1}Vol(supp\ h))^{\frac{2}{n}}\leq1\ and \ 2\rho_0\leq i_N,$$ \emph{where}
\[\rho_0=\left\{ \begin{array}{ll}
b^{-1}\sin^{-1}b(1-\alpha)^{-\frac{1}{n}}(\omega_n^{-1}Vol(supp\ h))^{\frac{1}{n}}\ \ \ &\ for\ b\ real,\\
(1-\alpha)^{-\frac{1}{n}}(\omega_n^{-1}Vol(supp\ h))^{\frac{1}{n}}\
\ \ &\ for\ b\ imaginary.
\end{array}
\right.\] \emph{ Here $\alpha$ is a free parameter, $0<\alpha<1$,
and}
$$C(n,\alpha)=\frac{1}{2}\pi\cdot2^{n-2}\alpha^{-1}(1-\alpha)^{-\frac{1}{n}}\frac{n}{n-1}\omega_n^{-\frac{1}{n}}.$$
\emph{For $b$ imaginary, we may omit the factor $\frac{1}{2}\pi$ in
the definition of $C(n,\alpha)$.}\\

The following lemma gives a proper form of the Sobolev inequality,
which can be found in \cite{12}. Here we outline the proof.\\
\\
\textbf{Lemma 2.3.} \emph{Let $M^n\subset N^{n+p}$ be a
$n(\geq3)$-dimensional closed submanifold in a Riemannian manifold
$N^{n+p}$ with codimension $p\geq1$. Denote by $i_N$ the positive
lower bound of the injective radius of $N$ restricted on $M$. Assume
$K_N\leq K_2$, where $K_2$ is a non-negative constant and let $f$ be
a non-negative $C^1$ function on $M$ satisfying
\begin{eqnarray}K_2(n+1)^{\frac{2}{n}}(\omega_n^{-1}Vol(supp\
f))^{\frac{2}{n}}\leq1,\end{eqnarray}
\begin{eqnarray}2K_2^{-\frac{1}{2}}\sin^{-1}K_2^{\frac{1}{2}}(n+1)^{\frac{1}{n}}(\omega_n^{-1}Vol(supp\
f))^{\frac{1}{n}}\leq i_N.\end{eqnarray} Then}
$$\parallel\nabla f\parallel^2_2\geq\frac{(n-2)^2}{4(n-1)^2(1+t)}\left[\frac{1}{C^2(n)}
\parallel f\parallel^2_{\frac{2n}{n-2}}-H^2_0\left(1+\frac{1}{t}\right)\parallel f\parallel^2_2\right],$$
\emph{where $H_0=\max_{x\in M}H$ and $C(n)=C(n,\frac{n}{n+1})$}\\
\\
\emph{Proof.}\ \ For all $g\in C^1(M)$, $g\geq0$ satisfying (1) and
(2), Lemma 2.2 implies \begin{eqnarray}\parallel
g\parallel_{\frac{n}{n-1}}\leq C(n)\int_M(\mid\nabla g\mid+Hg)d\mu.
\end{eqnarray}
Substituting $g=f^{\frac{2(n-1)}{n-2}}$ into $(3)$ gives
$$\left(\int_Mf^{\frac{2n}{n-2}}d\mu\right)^{\frac{n-1}{n}}\leq\frac{2(n-1)}{n-2}C(n)
\int_Mf^{\frac{n}{n-2}}\mid\nabla f\mid
d\mu+C(n)\int_MHf^{\frac{2(n-1)}{n-2}}d\mu.$$ By H\"{o}lder's
inequality, we get
$$\parallel f\parallel_{\frac{2n}{n-2}}\leq C(n)\left[\frac{2(n-1)}{n-2}\parallel\nabla f\parallel_2
+H_0\parallel f\parallel_2\right].$$ This implies
$$\parallel f\parallel^2_{\frac{2n}{n-2}}\leq C^2(n)\left[\frac{4(n-1)^2(1+t)}{(n-2)^2}\parallel\nabla f\parallel^2_2
+H^2_0\left(1+\frac{1}{t}\right)\parallel f\parallel^2_2\right],$$
which is desired.

\section{An estimate of the mean curvature by its $L^{n+2}$-norm}

\hspace*{6mm}In this section, we prove the following theorem, which
plays an
important role in the proof of Theorem 1.1.\\
\\
\textbf{Theorem 3.1.} \emph{ Suppose that $F_t:M^n\longrightarrow
N^{n+1}$ $(n\geq 3)$ is a mean curvature flow for $t\in[0,T_0]$, and
the second fundamental form $A$ is uniformly bounded on time
interval $[0,T_0]$. Then}
$$\max_{(x,t)\in M\times [\frac{T_0}{2},T_0]} H^2(x,t)\leq C_3\left(\int^{T_0}_0\int_{M_t}
 | H|^{n+2}d\mu dt\right)^{\frac{2}{n+2}},$$
\emph{where $C_3$ is a constant depending only on $n$, $T_0$,
$\sup_{(x,t)\in M\times[0,T_0]}|A|$, $K_1$, $K_2$ and the injectivity radius lower bound $i_N>0$ of $N$.}\\
\\
\emph{Proof.}\ \ The evolution equation of $H^2$ is
$$\frac{\partial}{\partial t}H^2=\triangle H^2-2\mid\nabla H\mid^2+2H^2\mid A\mid^2+2H^2\bar{R}ic(\nu,\nu).$$
Since $\mid A\mid$ is bounded, we obtain
\begin{eqnarray}\frac{\partial}{\partial t}H^2\leq\triangle H^2+\beta
H^2,
\end{eqnarray}
where $\beta$ is a positive constant depending only on $n$,
$\sup_{(x,t)\in M\times[0,T_0]}|A|$ and $K_2$. For $0<R<R'<\infty$
and $x\in M$, we set
\[\eta=\left\{ \begin{array}{ll}
1&\ \ \ \ \ \ x\in B_{g(0)}(x,R),\\
\eta\in [0,1]\ and \ |\nabla\eta|_{g(0)}\leq\frac{1}{R'-R}&\ \ \ \ \ \ x\in B_{g(0)}(x,R')\setminus B_{g(0)}(x,R),\\
0&\ \ \ \ \ \ x\in M\setminus B_{g(0)}(x,R').
\end{array}
\right.\] Since $supp\ \eta\subseteq B_{g(0)}(x,R')$, $\eta$
satisfies (1) and (2) with respect to $g(0)$ for $R'$ sufficiently
small. On the other hand, the area of some fixed subset in $M$ is
non-increasing along the mean curvature flow, hence $\eta$ satisfies
(1) and (2) with respect to each $g(t)$ for $t\in [0,T_0]$.

Fix $R'>0$ sufficiently small, for any point $x\in M_t$, we denote
by $B(R')=B_{g(0)}(x,R')$ the geodesic ball with radius $R'$
centered at $x$ with respect to the metric $g(0)$. Putting
$f=|H|^2$, then for any $p\geq 2$, the inequality $(4)$ implies
$$\frac{1}{p}\frac{\partial}{\partial t}\int_{B(R')} f^p\eta^2\leq\int_{B(R')}\eta^2 f^{p-1}\triangle f+\int_{B(R')}\beta f^p\eta^2
+\frac{1}{p}\int_{B(R')}f^p\eta^2\frac{\partial}{\partial
t}d\upsilon_t.$$ Integrating by parts yields
\begin{eqnarray*}
\int_{B(R')}\eta^2 f^{p-1}\triangle
f&=&-\frac{4(p-1)}{p^2}\int_{B(R')}\mid\nabla(f^{\frac{p}{2}}\eta)\mid^2+\frac{4}{p^2}\int_{B(R')}\mid\nabla\eta\mid^2f^p\\
&&+\frac{4(p-2)}{p^2}\int_{B(R')}\nabla(f^{\frac{p}{2}}\eta)f^{\frac{p}{2}}\nabla\eta\\
&\leq&-\frac{2}{p}\int_{B(R')}\mid\nabla(f^{\frac{p}{2}}\eta)\mid^2+\frac{2}{p}\int_{B(R')}\mid\nabla\eta\mid^2f^p.
\end{eqnarray*}
Thus
\begin{eqnarray*}
\frac{1}{p}\frac{\partial}{\partial t}\int_{B(R')}
f^p\eta^2&\leq&-\frac{2}{p}\parallel\nabla(f^{\frac{p}{2}}\eta)\parallel^2_2+\beta\parallel
f^{\frac{p}{2}}\eta\parallel^2_2\\
&&+\frac{2}{p}\int_{B(R')}\mid\nabla\eta\mid^2f^p+\frac{1}{p}\int_{B(R')}f^p\eta^2\frac{\partial}{\partial
t}d\upsilon_t\\
&\leq&-\frac{2}{p}\parallel\nabla(f^{\frac{p}{2}}\eta)\parallel^2_2+\beta\parallel
f^{\frac{p}{2}}\eta\parallel^2_2+\frac{2}{p}\int_{B(R')}\mid\nabla\eta\mid^2f^p.
\end{eqnarray*}
Hence \begin{eqnarray}\frac{\partial}{\partial t}\int_{B(R')}
f^p\eta^2+\int_{B(R')}\mid\nabla(f^{\frac{p}{2}}\eta)\mid^2\leq2\int_{B(R')}\mid\nabla\eta\mid^2f^p+\beta
p\int_{B(R')} f^p\eta^2.\end{eqnarray}

For any $0<\tau<\tau'<T_0$, define a function $\psi$ on $[0,T_0]$ by\\
\[\psi(t)=\left\{ \begin{array}{ll}
0&\ \ \ \ \ \ 0\leq t\leq \tau,\\
\frac{t-\tau}{\tau'-\tau}&\ \ \ \ \ \ \tau\leq t\leq \tau',\\
1&\ \ \ \ \ \ \tau'\leq t\leq T_0.
\end{array}
\right.\] Multiplying $(5)$ by $\psi(t)$ gives
\begin{eqnarray}&&\frac{\partial}{\partial t}\bigg(\psi\int_{B(R')}
f^p\eta^2\bigg)+\psi\int_{B(R')}\mid\nabla(f^{\frac{p}{2}}\eta)\mid^2\nonumber\\
&\leq&2\psi\int_{B(R')}\mid\nabla\eta\mid^2f^p+(\beta
p\psi+\psi')\int_{B(R')} f^p\eta^2.\end{eqnarray} By integrating
$(6)$ on $[\tau,t]$  we obtain
\begin{eqnarray*}&&\int_{B(R')}
f^p\eta^2+\int^t_{\tau'}\int_{B(R')}\mid\nabla(f^{\frac{p}{2}}\eta)
\mid^2\\
&\leq&2\int^{T_0}_\tau\int_{B(R')}\mid\nabla\eta\mid^2f^p
+\Big(\beta p+\frac{1}{\tau'-\tau}\Big)\int^{T_0}_\tau\int_{B(R')}
f^p\eta^2.\end{eqnarray*} For $R'$ sufficiently small, the following
Sobolev inequality holds:
\begin{eqnarray*}
\left(\int_{B(R')}
f^{\frac{pn}{n-2}}\eta^{\frac{2n}{n-2}}\right)^{\frac{n-2}{n}}&=&\parallel
f^{\frac{p}{2}}\eta\parallel^2_{\frac{2n}{n-2}}\\
&\leq&\frac{4(n-1)^2(1+s)C^2(n)}{(n-2)^2}\parallel\nabla(f^{\frac{p}{2}}\eta)\parallel^2_2\\
&&+H^2_0C^2(n)\Big(1+\frac{1}{s}\Big)\parallel
f^{\frac{p}{2}}\eta\parallel^2_2.
\end{eqnarray*}
Hence
\begin{eqnarray*}
&&\int^{T_0}_{\tau'}\int_{B(R')}
f^{p(1+\frac{2}{n})}\eta^{2+\frac{1}{n}}\nonumber\\
&\leq&\int^{T_0}_{\tau'}\left(\int_{B(R')}
f^{p}\eta^2\right)^{\frac{2}{n}}\left(\int_{B(R')} f^{\frac{np}{n-2}}\eta^{\frac{2n}{n-2}}\right)^{\frac{n-2}{n}} \\
&\leq&\max_{t\in [\tau',T_0]}\left(\int_{B(R')}
f^{p}\eta^2\right)^{\frac{2}{n}}\times\int_\tau^{T_0}\Big[\frac{4(n-1)^2(1+s)C^2(n)}{(n-2)^2}\parallel\nabla(f^{\frac{p}{2}}\eta)\parallel^2_2\\
&&+H^2_0C^2(n)\Big(1+\frac{1}{s}\Big)\parallel
f^{\frac{p}{2}}\eta\parallel^2_2\Big]\\
&\leq&C_4\max_{t\in [\tau',T_0]}\left(\int_{B(R')}
f^{p}\eta^2\right)^{\frac{2}{n}}\times\int_\tau^{T_0}\Big[\parallel\nabla(f^{\frac{p}{2}}\eta)\parallel^2_2+\parallel
f^{\frac{p}{2}}\eta\parallel^2_2\Big]\\&\leq&
C_4\left[2\int^{T_0}_{\tau}\int_{B(R')}\mid\nabla\eta\mid^2f^p+\Big(\beta
p+\frac{1}{\tau'-\tau}\Big)\int^{T_0}_\tau\int_{B(R')}
f^p\eta^2\right]^{1+\frac{2}{n}},\end{eqnarray*} where we put $s=1$
and $C_4$ is a constant depending on $n$ and $\sup_{(x,t)\in
M\times[0,T_0]}|A|$.

Note that $|\nabla\eta|_{g(t)}\leq|\nabla\eta|^2_{g(0)}e^{lt}$,
where $l=\max_{0\leq t\leq T_0}||\frac{\partial g}{\partial
t}||_{g(t)}$. Thus
$$\int^{T_0}_{\tau}\int_{B(R')}\mid\nabla\eta\mid^2f^p\leq
\int^{T_0}_{\tau}\int_{B(R')}\Big(\mid\nabla\eta\mid_{g(0)}e^{\frac{1}{2}lt}\Big)^2f^p\leq
\frac{e^{C_5T_0}}{(R'-R)^2}\int^{T_0}_{\tau}\int_{B(R')}f^p,$$ where
$C_5$ is a constant depending on $n$ and $\sup_{(x,t)\in
M\times[0,T_0]}|A|^2$. Then it follows that
\begin{eqnarray}
\int^{T_0}_{\tau}\int_{B(R)} f^{p(1+\frac{2}{n})}d\mu_tdt&\leq&
C_4\left(\beta p+\frac{1}{\tau'-\tau}+\frac{2e^{C_5T_0}}{(R'-R)^2}\right)^{1+\frac{2}{n}}\nonumber\\
&&\times\left(\int^{T_0}_\tau\int_{B(R')}
f^pd\mu_tdt\right)^{1+\frac{2}{n}}.
\end{eqnarray}
Putting $L(p,t,R)=\int^{T_0}_{t}\int_{B(R)}f^{p}$,  we obtain from
$(7)$
\begin{eqnarray}L\Big(p\Big(1+\frac{2}{n}\Big),\tau',R\Big)\leq
C_4\left(\beta
p+\frac{1}{\tau'-\tau}+\frac{2e^{C_5T_0}}{(R'-R)^2}\right)^{1+\frac{2}{n}}
L(p,\tau,R')^{1+\frac{2}{n}}.\end{eqnarray} We set
$\mu=1+\frac{2}{n}$, $p_k=\frac{n+2}{2}\mu^k$,
$\tau_k=\Big(1-\frac{1}{\mu^{k+1}}\Big)t$ and
$R_k=\frac{R}{2}\Big(1+\frac{1}{\mu^{k/2}}\Big)$, where
$k=0,1,2,\cdots$. Then it follows from $(8)$ that
\begin{eqnarray*}L(p_{k+1},\tau_{k+1},R_{k+1})^{\frac{1}{p_{k+1}}}&\leq&
C_4^{\frac{1}{p_{k+1}}}
\Big[\frac{(n+2)\beta}{2}+\frac{\mu^2}{\mu-1}\cdot\frac{1}{t}+\frac{4e^{C_5T_0}}{R'^2}\cdot\frac{\mu}
{(\sqrt{\mu}-1)^2}\Big]^{\frac{1}{p_k}}\\
&&\times\mu^{\frac{k}{p_{k}}}L(p_{k},\tau_{k},R_k)^{\frac{1}{p_{k}}}.\end{eqnarray*}
Hence for any $m\geq1$,
\begin{eqnarray*}
&&L(p_{m+1},\tau_{m+1},R_{m+1})^{\frac{1}{p_{m+1}}}\\
&\leq& C_4^{\sum^m_{k=0}\frac{1}{p_{k+1}}}
\Big[\frac{(n+2)\beta}{2}+\frac{\mu^2}{\mu-1}\cdot\frac{1}{t}+\frac{4e^{C_5T_0}}{R'^2}\cdot
\frac{\mu}{(\sqrt{\mu}-1)^2}\Big]^{\sum^m_{k=0}\frac{1}{p_{k}}}\\
&&\times\mu^{\sum^m_{k=0}\frac{k}{p_{k}}}L(p_{0},\tau_{0},R_0)^{\frac{1}{p_{0}}}.\end{eqnarray*}
As $m\rightarrow+\infty$, we conclude
\begin{eqnarray}f(x,t)\leq
C_6^{\frac{n}{n+2}}\Big(C_6+\frac{1}{t}+\frac{e^{C_5T_0}}{R'^2}\Big)\Big(1+\frac{2}{n}\Big)^{\frac{n}{2}}
\left(\int^{T_0}_0\int_{M_t}f^{\frac{{n+2}}{2}}\right)^{\frac{2}{n+2}},\end{eqnarray}
where $C_6$ is a positive constant depending on $n$, $\sup_{M\times
[0,T]}|A|$, $K_1$ and $K_2$.

Note that we choose $R'$ sufficient small such that
\begin{eqnarray}K_2(n+1)^{\frac{2}{n}}(\omega_n^{-1}Vol_{g(0)}(B(R'))^{\frac{2}{n}}\leq1\end{eqnarray}
and
\begin{eqnarray}2K_2^{-\frac{1}{2}}\sin^{-1}K_2^{\frac{1}{2}}(n+1)^{\frac{1}{n}}
(\omega_n^{-1}Vol_{g(0)}(B(R'))^{\frac{1}{n}}\leq i_N.\end{eqnarray}
For $g(0)$, there is a non-positive constant $K$ depending on $n$,
$\max_{x\in M_0}|A|$, $K_1$ and $K_2$ such that the sectional
curvature of $M_0$ is bounded from below by $K$. By the
Bishop-Gromov volume comparison theorem, we have
$$Vol_{g(0)}(B(R'))\leq Vol_K(B(R')),$$ where $Vol_K(B(R'))$ is the
volume of the ball with radius $R'$ in the $n$-dimensional complete
simply connected space form with constant curvature $K$. Let $R'$ be
the largest number satisfying
\begin{eqnarray}K_2(n+1)^{\frac{2}{n}}(\omega_n^{-1}Vol_K(B(R'))^{\frac{2}{n}}\leq1\nonumber\end{eqnarray}
and
\begin{eqnarray}2K_2^{-\frac{1}{2}}\sin^{-1}K_2^{\frac{1}{2}}(n+1)^{\frac{1}{n}}
(\omega_n^{-1}Vol_K(B(R'))^{\frac{1}{n}}\leq
i_N.\nonumber\end{eqnarray} Then $R'$ only depends on $n$, $K_1$,
$K_2$, $i_N$ and $\sup_{(x,t)\in M\times[0,T_0]}|A|$, and
$Vol_{g(0)}(B(R'))$ satisfies (10) and (11). This implies that
$$\max_{(x,t)\in M\times[\frac{T_0}{2},T_0]}H^2(x,t)\leq C_3\left(\int^{T_0}_0\int_{M_t}| H|^{n+2}d\mu dt\right)^{\frac{2}{n+2}},$$
where $C_3$ is a constant depending on $n$, $T_0$, $\sup_{(x,t)\in
M\times[0,T_0]}|A|$, $K_1$, $K_2$ and $i_N$, which is desired.\\

\section{Mean curvature flow with finite total mean curvature}

\hspace*{6mm}In this section we combine the above results to prove Theorem 1.1.\\
\\
\emph{Proof of Theorem 1.1.}\ \ It is sufficient to prove the
theorem for $\alpha=n+2$ since by H\"{o}lder's inequality, $||
H||_{\alpha,M\times [0,T)}<\infty$ implies $|| H||_{n+2,M\times
[0,T)}<\infty$ if $\alpha>n+2$. We argue by contradiction.

Suppose that the solution to the mean curvature flow can't be
extended over $T$, then $| A|$ becomes unbounded as $t\rightarrow
T$. Since $h_{ij}\geq -C$, we get $\sum_{i,j}(h_{ij}+C)^2\leq
C_7[tr(h_{ij}+C)]^2$, where $C_7$ is a constant depending only on
$n$. Since $| A|^2$ is unbounded, we have $\sum_{i,j}(h_{ij}+C)^2$
is unbounded. This together with
$$[tr(h_{ij}+C)]^2=(H+nC)^2=H^2+2nCH+n^2C^2$$
implies that $H^2$ is unbounded. Namely, $$\sup_{(x,t)\in M\times
[0,T)}H^2(x,t)=\infty.$$ Choose an increasing time sequence
$t^{(i)}$, $i=1,2,\cdots$, such that $\lim_{i\rightarrow
\infty}t^{(i)}=T$. We take a sequence of points $x^{(i)}\in M$
satisfying
$$H^2(x^{(i)},t^{(i)})=\max_{(x,t)\in M\times [0,t^{(i)}]} H^2(x,t).$$
Then $\lim_{i\rightarrow \infty} H^2(x^{(i)},t^{(i)})=\infty$.

Putting $Q^{(i)}=H^2(x^{(i)},t^{(i)})$, we have
$\lim_{i\rightarrow\infty}Q^{(i)}=\infty$. This together with
$\lim_{i\rightarrow \infty}t^{(i)}=T>0$ implies that there exists a
positive integer $i_0$ such that $Q^{(i)}t^{(i)}\geq1$ and
$Q^{(i)}\geq1$ for $i\geq i_0$. For $i\geq i_0$ and $t\in [0,1]$, we
consider the rescaled flows
$$F^{(i)}(t)=F\left(\frac{t-1}{Q^{(i)}}+t^{(i)}\right):(M,g^{(i)}(t))\longrightarrow(N,Q^{(i)}h).$$
Let $H_{(i)}$ and $g^{(i)}(t)=F^{(i)}(t)^*(Q^{(i)}h)$ be the mean
curvature of $F^{(i)}(t)$ and the induced metric on $M$ induced by
$F^{(i)}(t)$ respectively. Then $F^{(i)}(t):M\rightarrow
\mathbb{R}^{n+1}$ is still a solution to the mean curvature flow on
$t\in[0,1]$. Since $F_t$ satisfies $h_{ij}\geq-C$ for $(x,t)\in
M\times[0,T)$, we have
$$H^2_{(i)}(x,t)\leq 1 \ \ on\ \ M\times [0,1] ,$$
\begin{eqnarray}h^{(i)}_{jk}\geq -\frac{C}{\sqrt{Q^{(i)}}}\ \ on\ \
M\times [0,1],\end{eqnarray} where $A^{(i)}=h^{(i)}_{jk}$ is the
second fundamental form of $F^{(i)}(t)$. The inequality in $(12)$
gives $h^{(i)}_{jk}+\frac{C}{\sqrt{Q^{(i)}}}\geq 0$. Hence
$$h^{(i)}_{jk}+\frac{C}{\sqrt{Q^{(i)}}}\leq tr\left(h^{(i)}_{jk}+\frac{C}{\sqrt{Q^{(i)}}}\right)
\leq H_{(i)}+\frac{nC}{\sqrt{Q^{(i)}}},$$ which implies that
$h^{(i)}_{jk}\leq H_{(i)}+\frac{(n-1)C}{\sqrt{Q^{(i)}}}$. Since
$Q^{(i)}\geq1$ for $i\geq i_0$, it follows that $|A^{(i)}|\leq C_8$,
where $C_8$ is a constant independent of $i$ for $i\geq i_0$.

We consider the sequence $(M,g^{(i)}(t),x^{(i)})$, $t\in[0,1]$. It
follows from \cite{4} that there is a subsequence of
$(M^{(i)},g^{(i)}(t),x^{(i)})$ converges to a Riemannian manifold
$(\widetilde{M},\widetilde{g}(t),\widetilde{x})$, and the
corresponding subsequence of immersions $F^{(i)}(t)$ converges to an
immersion $\widetilde{F}(t):\widetilde{M}\rightarrow
\mathbb{R}^{n+1}$£¬ $t\in[0,1]$.

Since $(N,h)$ has bounded geometry and $Q^{(i)}\geq 1$ for $i\geq
i_0$, $(N,Q^{(i)}h)$ also has bounded geometry with the same
bounding constants as $(N,h)$ for each $i\geq i_0$. It follows from
Theorem 3.1 that
$$\max_{(x,t)\in M^{(i)}\times[\frac{1}{2},1]}H^2_{(i)}(x,t)\leq C_9\left(\int^{1}_0\int_{M_t}| H_{(i)}|^{n+2}d\mu_{g^{(i)}(t)}
dt\right)^{\frac{2}{n+2}},$$ where $C_9$ is a constant independent
of $i$  for $i\geq i_0$. Hence
\begin{eqnarray} \max_{(x,t)\in
\widetilde{M}\times[\frac{1}{2},1]}\widetilde{H}^2(x,t)&\leq
&\lim_{i\rightarrow \infty}C_9\left(\int^{1}_0\int_{M_t}|
H_{(i)}|^{n+2}d\mu_{g^{(i)}(t)} dt\right)^{\frac{2}{n+2}}\nonumber\\
&\leq&\lim_{i\rightarrow
\infty}C_9\left(\int^{t^{(i)}+(Q^{(i)})^{-1}}_{t^{(i)}}\int_{M_t}|
H|^{n+2}d\mu dt\right)^{\frac{2}{n+2}}.
\end{eqnarray}
The equality in $(13)$ holds because $\int^T_0\int_M  H^{n+2} d\mu
dt<+\infty$ and $\lim_{i\rightarrow \infty}(Q^{(i)})^{-1}=0$. On the
other hand, according to the choice of the points, we have
$$\widetilde{H}^2(\widetilde{x},1)=\lim_{i\rightarrow \infty}H^2_{(i)}(x^{(i)},1)=1.$$
This is a contradiction. We complete the proof of Theorem 1.1.\\

We can prove Theorem 1.2 by a similar method.\\\\
\emph{Proof of Theorem 1.2.} Since $H^2>n^2K_1$ is a strict
inequality and $M$ is compact, there is a positive constant
$\varepsilon$ such that $H^2\geq n^2K_1+\varepsilon$. From Lemma 3.1
we know that $H^2\geq n^2K_1+\varepsilon$ is preserved along the
flow. Moreover, $H^2>n^2K_1$ at $t=0$ implies $|A|^2\leq C_{10}
H^2$, where $C_{10}$ is a constant. We put $f_0=\frac{|A|^2}{H^2}$,
then we can obtain the evolution of $f_0$ by Lemma 5.2 in \cite{7}:
\begin{eqnarray*}
\frac{\partial}{\partial t}f_0&=&\triangle
f_0+\frac{2}{H}\langle\nabla_lH,\nabla_lf_0\rangle-\frac{2}{H^4}|\nabla_iHh_{kl}-\nabla_ih_{kl}H|^2\\
&&-\frac{1}{H^2}[4(h^{ij}h_{jl}{\bar{R}^l_{mi}\
^m}-h^{ij}h^{lm}\bar{R}_{iljm})+h^{ij}(\bar{\nabla}_j\bar{R}_{0li}\
^l+\bar{\nabla}_l\bar{R}_{0ij}\ ^l)].
\end{eqnarray*}
This implies that
$$\frac{\partial}{\partial t}f_0\leq\triangle
f_0+C_{11}f_0+C_{12},$$ where $C_{11}$ and $C_{12}$ are constants
independent of $t$. By the maximum principle and the finiteness of
$T$, there exists some positive constant $C_{13}$ independent of $t$
such that $|A|^2\leq C_{13} H^2$ for $t\in[0,T)$.

We only need to prove the theorem for $\alpha=n+2$, and we still
argue by contradiction. If the solution to the mean curvature flow
can't be extended over time $T$, then $|A|^2$ becomes unbounded as
$t\rightarrow T$, and $|A|^2\leq C_{13} H^2$ implies that $H^2$ also
becomes unbounded. Let $(x^{(i)},t^{(i)})$, $Q^{(i)}$, $F^{(i)}(t)$,
$g^{(i)}(t)$ and $(\widetilde{M},\tilde{g}(t),\tilde{x})$ be the
same as we denoted in the proof of Theorem 1.1, $A^{(i)}$ and
$H_{(i)}$ be the second fundamental form and mean curvature of the
immersion $F^{(i)}(t)$ respectively. Then $|A^{(i)}|^2\leq C_{13}
|H_{(i)}|^2$ for $(x,t)\in M\times [0,1]$, which implies that
$A^{(i)}$ is bounded by a constant independent of $i$, for $t\in
[0,1]$. Then by the conclusion of Theorem 3.1, we have
$$\max_{(x,t)\in M^{(i)}\times[\frac{1}{2},1]}H^2_{(i)}(x,t)\leq C_{14}
\left(\int^{1}_0\int_{M_t}| H_{(i)}|^{n+2}d\mu_{g^{(i)}(t)}
dt\right)^{\frac{2}{n+2}},$$ where $C_{14}$ is a constant
independent of $i$. Then by a similar process to the proof of
Theorem 1.1, we
can get a contradiction to complete the proof.\\
\\
\textbf{Remark 4.1.} By a similar argument, we can extend the mean
curvature flow in the case where $(M,g)$ is a complete non-compact
Riemannian manifold. But in that case, the condition $\alpha\geq
n+2$ have to be changed to
$\alpha=n+2$, since the H\"{o}lder's inequality doesn't hold.\\
\\
\textbf{Remark 4.2.} In \cite{17}, we have investigated the integral
conditions to extend mean curvature flow where $M^n$ is a
submanifold in $N^{n+p}$ with codimension $p\geq2$.

%-----------------------------------------------------------------------------------------------------------------

%-----------------------------------------------------------------------------------------------------------------

Center of Mathematical Sciences

Zhejiang University

Hangzhou 310027

China\\

E-mail address: xuhw@cms.zju.edu.cn; yf@cms.zju.edu.cn;
superzet@163.com

\end{document}